\documentclass[a4paper,12pt,reqno]{amsart}
\usepackage{graphics}
\newtheorem{theo}{Theorem}
\newtheorem{rem}{Remark}

\numberwithin{equation}{section}

\newcommand{\id}{\operatorname{id}}

\newcommand{\ad}{\operatorname{AD}}
\renewcommand{\r}{\operatorname{ROT}}
\newcommand{\h}{\operatorname{H}}

\begin{document}

\title[Alternating sign matrices]{Refined enumerations of alternating sign matrices: monotone $(d,m)$--trapezoids with prescribed top and bottom row}

\author[Ilse Fischer]{Ilse Fischer}

\thanks{Supported by the Austrian Science Foundation
    (FWF), grant number S9607--N13, in the framework of the National Research
    Network  ``Analytic Combinatorics and Probabilistic Number Theory''.}

\begin{abstract}
Monotone triangles are plane integer arrays of triangular shape
with certain monotonicity conditions along rows and diagonals.
Their significance is mainly due to the fact that they correspond
to $n \times n$ alternating sign matrices when prescribing
$(1,2,\ldots,n)$ as bottom row of the array. We define monotone
$(d,m)$--trapezoids as monotone triangles with $m$ rows where
the $d-1$ top rows are removed. (These
objects are also equivalent to certain partial alternating sign
matrices.) It is known that the number of
monotone triangles with bottom row $(k_1,\ldots,k_n)$ is given by
a polynomial $\alpha(n;k_1,\ldots,k_n)$ in the $k_i$'s. The main
purpose of this paper is to show that the number of monotone
$(d,m)$--trapezoids with prescribed top and bottom row appears as 
a coefficient in the expansion of a 
specialization of $\alpha(n;k_1,\ldots,k_n)$ with respect to a certain
polynomial basis. This settles a generalization of a
recent conjecture of Romik and the author. Among other things, the
result is used to express the number of monotone triangles with
bottom row $(1,2,\ldots,i-1,i+1,\ldots,j-1,j+1,\ldots,n)$ (which
is, by the standard bijection, also the number of $n \times n$
alternating sign matrices with given top two rows) in terms of the number of $n \times n$ alternating
sign matrices with prescribed top and bottom row, and, by a
formula of Stroganov for the latter numbers, to provide an
explicit formula for the first numbers. (A formula of this type
was first derived by Karklinsky and Romik using the relation of
alternating sign matrices to the six--vertex model.)
\end{abstract}

\maketitle

\section{Introduction}

A {\it monotone triangle} is a triangular array of integers of the following form
\begin{center}
\begin{tabular}{ccccccccccccccccc}
  &   &   &   &   &   &   &   & $a_{n,n}$ &   &   &   &   &   &   &   & \\
  &   &   &   &   &   &   & $a_{n-1,n-1}$ &   & $a_{n-1,n}$ &   &   &   &   &   &   & \\
  &   &   &   &   &   & $\dots$ &   & $\dots$ &   & $\dots$ &   &   &   &   &   & \\
  &   &   &   &   & $a_{3,3}$ &   & $\dots$ &   & $\dots$ &   & $a_{3,n}$ &   &   &   &   & \\
  &   &   &   & $a_{2,2}$ &   & $a_{2,3}$ &  &   $\dots$ &   & $\dots$   &  & $a_{2,n}$  &   &   &   & \\
  &   &   & $a_{1,1}$ &   & $a_{1,2}$ &   & $a_{1,3}$ &   & $\dots$ &   & $\dots$ &   & $a_{1,n}$ &   &   &
\end{tabular},
\end{center}
which is monotone increasing in northeast and in southeast
direction and strictly increasing along rows, that is $a_{i,j} \le
a_{i+1,j+1}$ for $1 \le i \le j < n$, $a_{i,j} \le a_{i-1,j}$ for
$1 < i \le j \le n$ and $a_{i,j} < a_{i,j+1}$ for $1 \le i \le j
\le n-1$. Monotone triangles with bottom row $(1,2,\ldots, n)$ are
said to be {\it complete} and correspond to $n \times n$
alternating sign matrices.  The
latter are defined as square matrices with entries in $\{0,1,-1\}$
such that in each row and column the sum of entries is $1$ and the non--zero entries alternate.
Next we display a complete
monotone triangle and its corresponding alternating sign matrix.
\begin{center}
\begin{tabular}{ccccccccccccccccc}
  &   &   &   &   &   &   &   & $4$ &   &   &   &   &   &   &   & \\
  &   &   &   &   &   &   & $3$ &   & $4$ &   &   &   &   &   &   & \\
  &   &   &   &   &   & $2$ &   & $3$ &   & $5$ &   &   &   &   &   & \\
  &   &   &   &   & $1$ &   & $2$ &   & $4$ &   & $6$ &   &   &   &   & \\
  &   &   &   & $1$ &   & $2$ &  &   $3$ &   & $5$   &  & $6$  &   &   &   & \\
  &   &   & $1$ &   & $2$ &   & $3$ &   & $4$ &   & $5$ &   & $6$ &   &   &
\end{tabular} \quad
$\begin{pmatrix} 0 & 0 & 0 & 1 & 0 & 0 \\
                 0 & 0 & 1 & 0 & 0 & 0 \\
                 0 & 1 & 0 & -1 & 1 & 0 \\
                 1 & 0 & -1 & 1 & -1 & 1 \\
                 0 & 0 & 1 & -1 & 1 &  0 \\
                 0 & 0 & 0 & 1 & 0 & 0
\end{pmatrix}$
\end{center}
In general, for a given complete monotone triangle $A=(a_{i,j})_{1 \le i \le j \le
n}$ with $n$ rows,  the corresponding $n \times n$
alternating sign matrix $M=(m_{i,j})_{1 \le i, j \le n}$ can be obtained as follows: we have $m_{i,j}=1$
if and only if
$$j \in
\{a_{n+1-i,n+1-i},a_{n+1-i,n+2-i},\ldots,a_{n+1-i,n}\} \setminus
\{a_{n+2-i,n+2-i},a_{n+2-i,n+3-i},\ldots,a_{n+2-i,n}\},
$$
$m_{i,j} = -1$ if and only if
$$j \in
\{a_{n+2-i,n+2-i},a_{n+2-i,n+3-i},\ldots,a_{n+2-i,n}\} \setminus
\{a_{n+1-i,n+1-i},a_{n+1-i,n+2-i},\ldots,a_{n+1-i,n}\}
$$
and $m_{i,j}=0$ else.

\medskip

The story of alternating sign matrices \cite{bressoud} started
back in the 1980s when Mills, Robbins and Rumsey \cite{mills}
conjectured that the number of $n \times n$ alternating sign
matrices is given by the following beautiful product formula
$$
\prod_{j=0}^{n-1} \frac{(3j+1)!}{(n+j)!}.
$$
It was more than $10$ years later when Zeilberger
\cite{zeilberger} finally succeeded in giving the first proof of
this formula.  Since then refined enumerations of alternating sign
matrices and enumerations of symmetry classes of alternating sign
matrices have been accomplished and fruitful relations to other
areas such as statistical physics (six--vertex model) and algebra
have been discovered and made the combinatorial analysis of these objects  to a major topic in
enumerative and algebraic combinatorics as well as in statistical physics.

\medskip

The central objects of this article are the following generalizations of
monotone triangles: for $1 \le d \le m$, a {\it
monotone $(d,m)$--trapezoid} is a monotone triangle with $m$ rows
where the first $d-1$ rows are removed. Thus,
\begin{center}
\begin{tabular}{ccccccccccccccccc}
  &   &   &   &   &   & $2$ &   & $3$ &   & $5$ &   &   &   &   &   & \\
  &   &   &   &   & $1$ &   & $2$ &   & $4$ &   & $6$ &   &   &   &   & \\
  &   &   &   & $1$ &   & $2$ &  &   $3$ &   & $5$   &  & $6$  &   &   &   &
\end{tabular}
\end{center}
is a monotone $(3,5)$--trapezoid. It is easy to translate these objects
into the language of alternating sign matrices: for $t \le n$, let
a $(t,n)$--partial alternating sign matrix be a $t \times n$
matrix with entries in $\{0,1,-1\}$ such that the non--zero
entries alternate in each row and column and the rowsums are
equal to $1$. Suppose $1 \le i_1 < i_2 < \ldots < i_d \le n$ and $1 \le s_1 < s_2 < \ldots < s_c \le n$. Then $(d,n-c)$--trapezoids with top
row $i_1,\ldots,i_d$ and whose bottom row consists of the numbers in
$\{1,2,\ldots,n\} \setminus \{s_{1},s_{2},\ldots,s_{c}\}$,
arranged in increasing order, are equivalent to
$(n-c-d,n)$--partial alternating sign matrices with the following
properties: the $j$--th columnsum is $1$ iff $j \notin
\{i_1,i_2,\ldots,i_d,s_1,s_2,\ldots,s_c\}$ and it is $-1$ iff $j
\in \{i_1,\ldots,i_d\} \cap \{s_1,\ldots s_c\}$. Otherwise, the
$j$--th columnsum is $0$ and the first non--zero entry (if there
is any) of this column is a $-1$ iff $j \in \{i_1,\ldots, i_d\}
\setminus \{ s_1, \ldots, s_c\}$. The $(2,6)$--partial
alternating sign matrix that corresponds to the monotone
$(3,5)$--trapezoid given above is
$$\begin{pmatrix}
                 1 & 0 & -1 & 1 & -1 & 1 \\
                 0 & 0 & 1 & -1 & 1 &  0
\end{pmatrix}\footnote{As the monotone $(3,5)$--trapezoid consists of rows $3-5$ of our introductory
example of a complete monotone triangle, the $(2,6)$--partial alternating sign matrix
is the fourth and fifth row of the corresponding alternating
sign matrix.}.$$

\medskip

Let $\alpha(n;k_1,\ldots,k_n)$ denote the number of monotone
triangles with bottom row $(k_1,\ldots,k_n)$. In \cite{fischergtriangle} it
was shown that
\begin{equation}
\label{operator} \alpha(n;k_1,\ldots,k_n) = \left( \prod_{1 \le p
< q \le n} (\id + E_{k_p} E_{k_q} - E_{k_q} ) \right) \prod_{1 \le
i < j \le n} \frac{k_j - k_i}{j -i},
\end{equation}
where $E_{x}$ denotes the shift operator, defined as $E_{x} p(x) =
p(x+1)$. Moreover note that the product of operators in this
formula is understood as the composition.
In particular, this operator formula implies that
$\alpha(n;k_{1},\ldots,k_{n})$ is a polynomial in
$(k_{1},\ldots,k_{n})$ of degree no greater than $n-1$ in each
$k_{i}$. For non--negative integers $c,d$ with $c+d \le n$, we
consider the following expansion.
\begin{multline}
\label{expansion}
\alpha(n;k_{1},\ldots,k_{c},c+1,c+2,\ldots,n-d,k_{n-d+1},k_{n-d+2},\ldots,k_{n})
\\
= \sum_{s_{1}=1}^{n} \sum_{s_{2}=1}^{n} \ldots \sum_{s_{c}=1}^{n}
\sum_{i_{1}=1}^{n} \sum_{i_{2}=1}^{n} \ldots \sum_{i_{d}=1}^{n}
A(n;s_{1},s_{2},\ldots,s_{c};i_{1},\ldots,i_{d}) \\
\times (-1)^{s_{1}+s_{2}+\ldots+s_{c}+c}
\binom{k_{1}-c-1}{s_{c}-1} \binom{k_{2}-c-1}{s_{c-1}-1} \cdots
\binom{k_{c}-c-1}{s_{1}-1} \\ \times
\binom{k_{n-d+1}-n+d-2+i_{1}}{i_{1}-1} \binom{k_{n-d+2}-n+d-2+i_{2}}{i_{2}-1}
\cdots \binom{k_{n}-n+d-2+i_{d}}{i_{d}-1}.
\end{multline}
The main subject of this paper is the following generalization of
Conjecture~7 from \cite{fischerromik}, which provides a combinatoral interpretation of the 
coefficients $A(n;s_{1},s_{2},\ldots,s_{c};i_{1},\ldots,i_{d})$. (In fact, this theorem also implies
Conjecture~5 from \cite{fischerromik}, see Section~\ref{remarks}.)

\begin{theo}
\label{7} For $1 \le s_{1} < s_{2} < \cdots < s_{c} \le n$ and $1 \le i_{1} < i_{2}
< \cdots < i_{d} \le n$ the coefficient
$A(n;s_{1},\ldots,s_{c};i_{1},\ldots,i_{d})$ is the number of
monotone $(d,n-c)$--trapezoids with $(i_{1},\ldots,i_{d})$ as top
row and whose bottom row consists of the numbers in
$\{1,2,\ldots,n\} \setminus \{s_{1},s_{2},\ldots,s_{c}\}$,
arranged in increasing order.
\end{theo}

Its proof is the content of
Sections~\ref{rotreflect}--\ref{complet}.
In Section~\ref{relate}, we use the identity
\begin{equation*}
 \alpha(n;k_1,\ldots,k_n) = (-1)^{n-1}
\alpha(n;k_2,\ldots,k_n,k_1-n)
\end{equation*}
from \cite[Lemma~5]{newproof} to express the coefficient
$$A(n;s_1,\ldots,s_c;i_1,\ldots,i_{d})$$ in terms of the coefficients
$$A(n;s_1,\ldots,s_{c-t};i_1,\ldots,i_{d+t})$$ where $t$ is a fixed non--negative integer no greater than $c$ and
$s_{c-l+1} \le i_{d+l} \le n$ for $1 \le l \le t$. In the special case
$c=2$, $d=0$ and $t=1$ this relates
$A(n;s_1,s_2;-)=:A_{n,s_{1},s_{2}}$ to
$A(n;s_1;i_1)=:B_{n,s_{1},i_{1}}$. By the standard bijection
between alternating sign matrices and complete monotone triangles,
$B_{n,i,j}$ is the number of $n \times n$ alternating sign matrices
where the unique $1$ in the bottom row is in column $i$ and the
unique $1$ in the top row is in column $j$. Stroganov
\cite{stroganov} has derived a formula for this
doubly refined enumeration of alternating sign matrices. On the
other hand, the numbers $A_{n,i,j}$ have been studied in
\cite{fischerromik} and our relation to Stroganov's doubly refined
enumeration numbers finally enables us  to provide the
following formula for the numbers $A_{n,i,j}$:
\begin{multline}
\label{anij}
 A_{n,i,j} = \frac{1}{A_{n-1}} \sum_{i=1}^{l} \sum_{k=l-i+j}^{l-i+n}   (-1)^{n+i+k+l}
\binom{2n-2-j}{k-l+i-j} \\ \times \left( A_{n-1,l-1} (A_{n,k} -
A_{n,k-1}) + A_{n-1,k-1} (A_{n,l} -  A_{n,l-1}) \right),
\end{multline}
where
$$
A_{n,k} = \left\{ \begin{array}{cc} \binom{n+k-2}{n-1}
\frac{(2n-k-1)!}{(n-k)!} \prod\limits_{j=0}^{n-2}
\frac{(3j+1)!}{(n+j)!}
& 1 \le k \le n \\
0 & \text{otherwise} \end{array} \right.
$$
is the number of $n \times n$ alternating sign matrices that have
a $1$ in the $k$--th column of the first row and $A_{n}$ denotes
the total number of $n \times n$ alternating sign matrices. A
formula similar to \eqref{anij} was first derived by Karklinsky
and Romik \cite{2topbottom} by using the six--vertex model approach
to alternating sign matrices. Another, complicated, formula for
$A_{n,i,j}$ was conjectured in \cite[Conjecture~4]{fischerromik}.)

\medskip

In fact, the numbers $A_{n,i,j}$ also correspond to a certain
doubly refined enumeration of $n \times n$ alternating sign
matrices with respect to the two bottom (or, equivalently, two top) rows if $i < j$: let $k$ be a fixed integer with $i \le k \le j$.
Then $A_{n,i,j}$ is the number of $n \times n$ alternating sign
matrices such that $e_{i}+e_{j}-e_k$ is the $(n-1)$--st row and
$e_k$ is the bottom row where $e_p \in \mathbb{R}^n$ with
$(e_p)_q = \delta_{p,q}$. Alternatively and as explained above, it
is also the number of $(n-2,n)$--partial alternating sign matrices
such that the last non--zero entries (if there exist any) in
columns $i$ and $j$ are $-1$s.

\section{How does rotating and reflecting alternating sign matrices translate into the language of complete
  monotone triangles?}
\label{rotreflect}

The set of $n \times n$ alternating sign matrices is invariant
under the rotations of $90^{\circ}$ and under the reflection along
any of the four symmetry axes of the square. Here, we investigate
how these rotations and reflections translate into the language of
complete monotone triangles.

\medskip

Let $A=(a_{i,j})_{1 \le i \le j \le n}$ be a complete monotone
triangle with $n$ rows and $M$ be the corresponding $n \times n$
alternating sign matrix. We fix the following notation: the
sequence $(a_{l,l},a_{l-1,l},\ldots,a_{1,l})$ is said to be the
$l$--th SE--diagonal and the sequence
$(a_{1,l},a_{2,l+1},\ldots,a_{n-l+1,n})$ is said to be the $l$--th
NE--diagonal of the complete monotone triangle $A$.

\medskip

It is not hard to convince oneself that the complete monotone triangle
$B=(b_{i,j})_{1 \le i \le j \le n}$ which corresponds to the
alternating sign matrix that we obtain by reflecting $M$ along the
antidiagonal  is given by
$$
b_{i,j} = \text{$\#$ of elements $x$ in the $j$--th SE--diagonal
of $A$ with $x \ge i$}.
$$
We set $B=\ad(A)$. For the example in the introduction we obtain
the following monotone triangle and alternating sign
matrix.
\begin{center}
\begin{tabular}{ccccccccccccccccc}
  &   &   &   &   &   &   &   & $3$ &   &   &   &   &   &   &   & \\
  &   &   &   &   &   &   & $2$ &   & $4$ &   &   &   &   &   &   & \\
  &   &   &   &   &   & $1$ &   & $3$ &   & $6$ &   &   &   &   &   & \\
  &   &   &   &   & $1$ &   & $2$ &   & $5$ &   & $6$ &   &   &   &   & \\
  &   &   &   & $1$ &   & $2$ &  &   $4$ &   & $5$   &  & $6$  &   &   &   & \\
  &   &   & $1$ &   & $2$ &   & $3$ &   & $4$ &   & $5$ &   & $6$ &   &   &
\end{tabular}
 \quad
$\begin{pmatrix} 0 & 0 & 1 & 0 & 0 & 0 \\
                 0 & 1 & -1 & 1 & 0 & 0 \\
                 1 & -1 & 1 & -1 & 0 & 1 \\
                 0 & 1 & -1 & 0 & 1 & 0 \\
                 0 & 0 & 0 & 1 & 0 &  0 \\
                 0 & 0 & 1 & 0 & 0 & 0
\end{pmatrix}$
\end{center}

\medskip

Similarly, the complete monotone triangle $C=(c_{i,j})_{1 \le i \le j \le
n}$ which corresponds to the alternating sign matrix that we
obtain by rotating $M$ clockwise by $90^{\circ}$ is given by
$$
c_{i,j} = \text{$\#$ of elements $x$  in the $(n+1-j)$--th
NE--diagonal of $A$ with $x \le n+1-i$}.
$$
We set $C=\r(A)$. In the example this gives the following two
objects.
\begin{center}
\begin{tabular}{ccccccccccccccccc}
  &   &   &   &   &   &   &   & $3$ &   &   &   &   &   &   &   & \\
  &   &   &   &   &   &   & $3$ &   & $4$ &   &   &   &   &   &   & \\
  &   &   &   &   &   & $2$ &   & $4$ &   & $5$ &   &   &   &   &   & \\
  &   &   &   &   & $1$ &   & $3$ &   & $5$ &   & $6$ &   &   &   &   & \\
  &   &   &   & $1$ &   & $2$ &  &   $4$ &   & $5$   &  & $6$  &   &   &   & \\
  &   &   & $1$ &   & $2$ &   & $3$ &   & $4$ &   & $5$ &   & $6$ &   &   &
\end{tabular} \quad
$\begin{pmatrix} 0 & 0 & 1 & 0 & 0 & 0 \\
                 0 & 0 & 0 & 1 & 0 & 0 \\
                 0 & 1 & -1 & 0 & 1 & 0 \\
                 1 & -1 & 1 & -1 & 0 & 1 \\
                 0 & 1 & -1 & 1 & 0 &  0 \\
                 0 & 0 & 1 & 0 & 0 & 0
\end{pmatrix}$
\end{center}

\medskip

Finally, the complete monotone triangle $D=(d_{i,j})_{1 \le i \le j \le n}$
that corresponds to the alternating sign matrix that we obtain by
reflecting $M$ along the horizontal symmetry axis is uniquely
determined by
$$\{ a_{i,i},
a_{i,i+1},\ldots,a_{i,n},d_{n+2-i,n+2-i},d_{n+2-i,n+3-i},\ldots,d_{n+2-i,n}
\} = \{1,2,\ldots,n\}$$ for all $i \in \{1,2,\ldots,n\}$. Here we
set $D=\h(A)$. In this case the running example changes to
\begin{center}
\begin{tabular}{ccccccccccccccccc}
  &   &   &   &   &   &   &   & $4$ &   &   &   &   &   &   &   & \\
  &   &   &   &   &   &   & $3$ &   & $5$ &   &   &   &   &   &   & \\
  &   &   &   &   &   & $1$ &   & $4$ &   & $6$ &   &   &   &   &   & \\
  &   &   &   &   & $1$ &   & $2$ &   & $5$ &   & $6$ &   &   &   &   & \\
  &   &   &   & $1$ &   & $2$ &  &   $3$ &   & $5$   &  & $6$  &   &   &   & \\
  &   &   & $1$ &   & $2$ &   & $3$ &   & $4$ &   & $5$ &   & $6$ &   &   &
\end{tabular}\quad
$\begin{pmatrix} 0 & 0 & 0 & 1 & 0 & 0 \\
                 0 & 0 & 1 & -1 & 1 & 0 \\
                 1 & 0 & -1 & 1 & -1 & 1 \\
                 0 & 1 & 0 & -1 & 1 & 0 \\
                 0 & 0 & 1 & 0 & 0 &  0 \\
                 0 & 0 & 0 & 1 & 0 & 0
\end{pmatrix}$.
\end{center}

\medskip

Since reflecting an alternating sign matrix along the antidiagonal
is equivalent to first rotating it clockwise by $90^\circ$ and
then reflecting it along the horizontal symmetry axis we have
\begin{equation}
\label{identity} \ad = \h \circ \r.
\end{equation}

\medskip

However, in the following we do not need the interpretations of
the mappings $\ad, \r$ and $\h$ in terms of alternating sign
matrices. For our purpose, it suffices to show that these three
mappings are permutations of the set of complete monotone
triangles with $n$ rows (this is easy and left to the reader) and
that they satisfy \eqref{identity}. The latter is equivalent to
showing that for every complete monotone triangle $A=(a_{i,j})_{1
\le i \le j \le n}$ and $i \in \{1,2,\ldots,n\}$ the union of
$$
\{ \text{$\#$ of elements $x$ in $j$--th SE--diagonal of $A$ with
$x \ge i$} | i \le j \le n \}
$$
and $$\{ \text{$\#$ of elements $x$ in the $j$-th NE--diagonal of
$A$ with $x \le i-1$} | 1 \le j \le i-1\}$$ is equal to
$\{1,2,\ldots,n\}$. This follows from the following fact which can
be shown by induction with respect to $n$, see
Figure~\ref{kernel}: suppose $(p_1,p_2,\ldots,p_{i-1})$ is a
strictly decreasing sequence of positive integers and for $1 \le j
\le i-1$ delete the first $p_j$ entries from the $j$--th
NE--diagonal of a monotone triangle with $n$ rows. For $i \le j
\le n$, let $q_j$ be the length of the $j$--th SE--diagonal in the
remaining (partial) monotone triangle. Then
$$
\{p_1,\ldots,p_{i-1},q_i,q_{i+1},\ldots,q_n\} = \{1,2,\ldots,n\}.
$$

\begin{figure}[h]
\begin{center}
\begin{tabular}{ccccccccccccccccc}
  &   &   &   &   &   &   &   & $\circ$ &   &   &   &   &   &   &   & \\
  &   &   &   &   &   &   & $\bullet$ &   & $\circ$ &   &   &   &   &   &   & \\
  &   &   &   &   &   & $\bullet$ &   & $\bullet$ &   & $\circ$ &   &   &   &   &   & \\
  &   &   &   &   & $\bullet$ &   & $\bullet$ &   & $\circ$ &   & $\circ$ &   &   &   &   & \\
  &   &   &   & $\bullet$ &   & $\bullet$ &  &   $\bullet$ &   & $\circ$   &  & $\circ$  &   &   &   & \\
  &   &   & $\bullet$ &   & $\bullet$ &   & $\bullet$ &   & $\bullet$ &   & $\circ$ &   & $\circ$ &   &   &
\end{tabular}
\end{center}
\caption{$p_1=5$, $p_2=4$, $p_3=2$, $p_4=1$, $q_5=3$, $q_6=6$}
\label{kernel}
\end{figure}

\section{From monotone $(d,n-c)$--trapezoids with prescribed top and bottom row to
  other partial monotone triangles}

In this section we fix the two sequences $(i_1,\ldots,i_d)$ and $(s_1,\ldots,s_c)$  with 
$1 \le i_1 < i_2 < \ldots < i_d \le n$ and $1 \le s_1 < s_2 < \ldots < s_c \le n$.
We will use the mappings from the previous section to show
that the number of monotone $(d,n-c)$--trapezoids with
$(i_{1},\ldots,i_{d})$ as top row and whose bottom row consists of
the numbers in $\{1,2,\ldots,n\} \setminus
\{s_{1},s_{2},\ldots,s_{c}\}$, arranged in increasing order, is
equal to the number of (partial) monotone triangles with $n$ rows,
entries in $\{c+1,c+2,\ldots,n-d\}$ where for $1 \le l \le c$ the
first $s_{c+1-l}$ entries of the $l$--th NE--diagonal are missing
and for $n-d+1 \le l \le n$ the last $i_{l-n+d}$ entries of the
$l$--th SE--diagonal are missing.

\medskip

Indeed, fix such a monotone $(d,n-c)$--trapezoid and add entries
arbitrarily to make it up to a complete monotone triangle $A$ with
$n$ rows. As an example, we take the $(3,5)$--trapezoid on the
left and the completion on the right, i.e. $c=1$ and $n=6$.
\begin{center}
\begin{tabular}{ccccccccccccc}
  &   & $\bf 1$ &   & $\bf 3$ &   & $\bf 6$ &   &  \\
  & $\bf 1$ &   & $\bf 2$ &   & $\bf 5$ &   & $\bf 6$ &   \\
   $\bf 1$ &   & $\bf 2$ &  &   $\bf 4$ &   & $\bf 5$   &  & $\bf 6$
\end{tabular} \quad
\begin{tabular}{ccccccccccccccccc}
  &   &   &   &   &   &   &   &  $3$ &   &   &   &   &   &   &   & \\
  &   &   &   &   &   &   & $2$ &   & $4$ &   &   &   &   &   &   & \\
  &   &   &   &   &   & $\bf 1$ &   & $\bf 3$ &   & $\bf 6$ &   &   &   &   &   & \\
  &   &   &   &   & $\bf 1$ &   & $\bf 2$ &   & $\bf 5$ &   & $\bf 6$ &   &   &   &   & \\
  &   &   &   & $\bf 1$ &   & $\bf 2$ &  &   $\bf 4$ &   & $\bf 5$   &  & $\bf 6$  &   &   &   & \\
  &   &   & $1$ &   & $2$ &   & $3$ &   & $4$ &   & $5$ &   & $6$ &   &   &
\end{tabular}
\end{center}
Consider $B=\ad^{-1}(A)$. In our example $B$ is the following complete 
monotone triangle.
\begin{center}
\begin{tabular}{ccccccccccccccccc}
  &   &   &   &   &   &   &   & $4$ &   &   &   &   &   &   &   & \\
  &   &   &   &   &   &   & $\bf 3$ &   & $4$ &   &   &   &   &   &   & \\
  &   &   &   &   &   & $\bf 2$ &   & $\bf 3$ &   & $5$ &   &   &   &   &   & \\
  &   &   &   &   & $1$ &   & $\bf 2$ &   & $4$ &   & $6$ &   &   &   &   & \\
  &   &   &   & $1$ &   & $\bf 2$ &  &   $\bf 3$ &   & $5$   &  & $6$  &   &   &   & \\
  &   &   & $1$ &   & $\bf 2$ &   & $\bf 3$ &   & $4$ &   & $5$ &   & $6$ &   &   &
\end{tabular}
\end{center}
(The entries $c+1, c+2, \ldots,n-d$ are displayed in boldface as
they correspond to the original monotone $(3,5)$--trapezoid.) By the definition of 
$\ad$, the
elements greater than $n-d$ in $B$ are just the last $i_{l-n+d}$
elements in the $l$--th SE--diagonals where $n-d+1 \le l \le n$ and
the distribution of these entries is exactly determined by the
completion of the monotone $(d,n-c)$--trapezoid towards the top.
On the other hand, by \eqref{identity}, $B=\r^{-1} \circ \h^{-1}
(A)$. Thus, the $c$--th row of $\h^{-1} (A)$ (counted from the
top) is $(s_1,\ldots,s_c)$ and the first $c-1$ rows of this complete
monotone triangle correspond to the completion of the original monotone
$(d,n-c)$--trapezoid towards the bottom. Consequently (and by the
definition of $\r$), the elements no greater than $c$ of $B$ are
just the first $s_{c+1-l}$ elements in the $l$--th NE--diagonals
for $1 \le l \le c$ and the distribution of these entries is
determined by the completion of the monotone $(d,n-c)$--trapezoid
towards the bottom.

\medskip

If we delete in our example the first $s_{c+1-l}$ entries in the $l$--th
NE--diagonal for $1 \le l \le c$ and the last $i_{l-n+d}$ elements in
the $l$--th SE--diagonal for $n-d+1 \le l \le n$ (i.e. exactly the part of
the complete monotone triangles that corresponds to the completion of the original
monotone $(d,n-c)$--trapezoid towards top and bottom), we clearly obtain the
following partial monotone triangle.
\begin{center}
\begin{tabular}{ccccccccccccccccc}
  &   &   &   &   &   &   &   & \phantom{$4$} &   &   &   &   &   &   &   & \\
  &   &   &   &   &   &   & $\bf 3$ &   &  \phantom{$4$} &   &   &   &   &   &   & \\
  &   &   &   &   &   & $\bf 2$ &   & $\bf 3$ &   &  \phantom{$5$} &   &   &   &   &   & \\
  &   &   &   &   &  \phantom{$1$} &   & $\bf 2$ &   &  \phantom{$4$} &   &  \phantom{$6$} &   &   &   &   & \\
  &   &   &   &  \phantom{$1$} &   & $\bf 2$ &  &   $\bf 3$ &   &  \phantom{$5$}   &  &  \phantom{$6$}  &   &   &   & \\
  &   &   &  \phantom{$1$} &   & $\bf 2$ &   & $\bf 3$ &   &  \phantom{$4$} &   &  \phantom{$5$} &   &  \phantom{$6$} &   &   &
\end{tabular}
\end{center}


\medskip

For technical reasons (which become apparent later) we add the
entry $c+1$ at the beginning of the truncated $l$--th NE--diagonal
for $1 \le l \le c$ and add the entry $n-d$ at the end of the
truncated $l$--th SE--diagonal for $n-d+1 \le l \le n$. In our
example, we obtain
\begin{center}
\begin{tabular}{ccccccccccccccccc}
  &   &   &   &   &   &   &   & $3$ &   &   &   &   &   &   &   & \\
  &   &   &   &   &   &   & $\bf 3$ &   &  \phantom{$4$} &   &   &   &   &   &   & \\
  &   &   &   &   &   & $\bf 2$ &   & $\bf 3$ &   &  \phantom{$5$} &   &   &   &   &   & \\
  &   &   &   &   &  $2$ &   & $\bf 2$ &   &  $3$ &   &  \phantom{$6$} &   &   &   &   & \\
  &   &   &   &  \phantom{$1$} &   & $\bf 2$ &  &   $\bf 3$ &   &  \phantom{$5$}   &  &  \phantom{$6$}  &   &   &   & \\
  &   &   &  \phantom{$1$} &   & $\bf 2$ &   & $\bf 3$ &   &  $3$ &   &  \phantom{$5$} &   &  \phantom{$6$} &   &   &
\end{tabular}.
\end{center}
 This will certainly not destroy the monotonicity along diagonals, but may
 disturb the strict increase along rows. All in all we see that the number of
monotone $(d,n-c)$--trapezoids with fixed top and bottom row as
given above is the number of partial monotone triangles with $n$
rows and entries in $\{c+1,c+2,\ldots,n-d\}$, where for $1 \le l
\le c$ the first $s_{c+1-l}-1$ entries of the $l$--th NE--diagonal
are missing and the first entry of the truncated diagonal is equal
to $c+1$ and does not necessarily have to be strictly smaller than
its neighbour to the right and for $n-d+1 \le l \le n$ the last
$i_{l-n+d}-1$ entries in the $l$--th SE--diagonal are missing and
the last entry of the truncated diagonal is equal to $n-d$ and
does not necessarily have to be strictly greater than its
neighbour to the left.

\section{Completion of the proof of Theorem~\ref{7}}
\label{complet}

Suppose $A(l_1,\ldots,l_{n-1})$ is a function in
$(l_1,\ldots,l_{n-1})$ and $(k_1,\ldots,k_{n}) \in \mathbb{Z}^n$.
We define a summation operator
$$
\sum_{(l_1,\ldots,l_{n-1})}^{(k_1,\ldots,k_n)}
A(l_1,\ldots,l_{n-1})
$$
by induction with respect to $n$. If $n=0$ then the application of
the operator gives zero, for $n=1$ we set $\sum\limits_{(-)}^{k_1}
A = A$. If $n \ge 2$ then we define
\begin{multline*}
\sum_{(l_1,\ldots, l_{n-1})}^{(k_1,\ldots,k_n)} A(l_1,\ldots,l_{n-1}) =
\sum_{(l_1,\ldots,l_{n-2})}^{(k_1,\ldots,k_{n-1})} \sum_{l_{n-1} =
k_{n-1}}^{k_n} A(l_1,\ldots,l_{n-1}) \\ -
\sum_{(l_1,\ldots,l_{n-3})}^{(k_1,\ldots,k_{n-2})}
A(l_1,\ldots,l_{n-3},k_{n-1},k_{n-1}).
\end{multline*}
With this, we
obviously have
$$
\alpha(n;k_1,\ldots,k_n) =
\sum_{(l_1,\ldots,l_{n-1})}^{(k_1,\ldots,k_n)}
\alpha(n-1;l_1,\ldots,l_{n-1})
$$
for $n \ge 2$.
In fact, if we extend the definition of monotone triangles to
triangular integer arrays $(a_{i,j})_{1 \le i \le j \le n}$ with
weak increase along NE--diagonals and SE--diagonals and strict
increase along rows with the possible exception of the bottom row
which may only be weakly increasing then
$\alpha(n;k_1,\ldots,k_n)$ is also the number of these (extended) monotone
triangles with bottom row $(k_1,\ldots,k_n)$. Let $\delta_x = \id
- E^{-1}_x$ and set
$$
B(k_{1},\ldots,k_{n}) = \sum_{(l_1,\ldots,l_{n-1})}^{(k_1,\ldots,k_n)}
A(l_1,\ldots,l_{n-1}).
$$
 The recursion shows that
$$\delta_{k_n} B(k_1,\ldots,k_n) =
\sum_{(l_1,\ldots,l_{n-2})}^{(k_1,\ldots,k_{n-1})}
A(l_1,\ldots,l_{n-2},k_n)$$ and, more general for $d \ge 1$,
$$
\delta_{k_{n-d+1}} \delta_{k_{n-d+2}} \dots   \delta_{k_{n}} B(k_1,\ldots,k_n) =
\sum_{(l_1,\ldots,l_{n-d-1})}^{(k_1,\ldots,k_{n-d})}
A(l_1,\ldots,l_{n-d-1},k_{n-d+1},\ldots,k_{n}).
$$
Similarly, with $\Delta_x = E_x - \id$, it is not hard to see that
$$(-1)^{c} \Delta_{k_1} \Delta_{k_{2}} \dots \Delta_{k_{c}}
B(k_1,\ldots,k_n) =
\sum_{(l_{c+1},\ldots,l_{n-1})}^{(k_{c+1},\ldots,k_{n})}
A(k_{1},\ldots,k_{c},l_{c+1},\ldots,l_{n-1}).
$$
If we combine these two facts we see that
\begin{multline}
\label{delta}
 \Delta_{k_1} \Delta_{k_{2}} \dots \Delta_{k_{c}} \delta_{k_{n-d+1}}
 \delta_{k_{n-d+2}} \dots   \delta_{k_{n}} (-1)^{c}
 B(k_{1},\ldots,k_{n}) \\
= \sum_{(l_{c+1},\ldots,l_{n-d-1})}^{(k_{c+1},\ldots,k_{n-d})} A(k_{1},\ldots,k_{c},l_{c+1},\ldots,l_{n-d-1},k_{n-d+1},\ldots,k_{n}).
\end{multline}

\medskip

\begin{figure}
\begin{center}
\mbox{\scalebox{0.50}{%
    \includegraphics{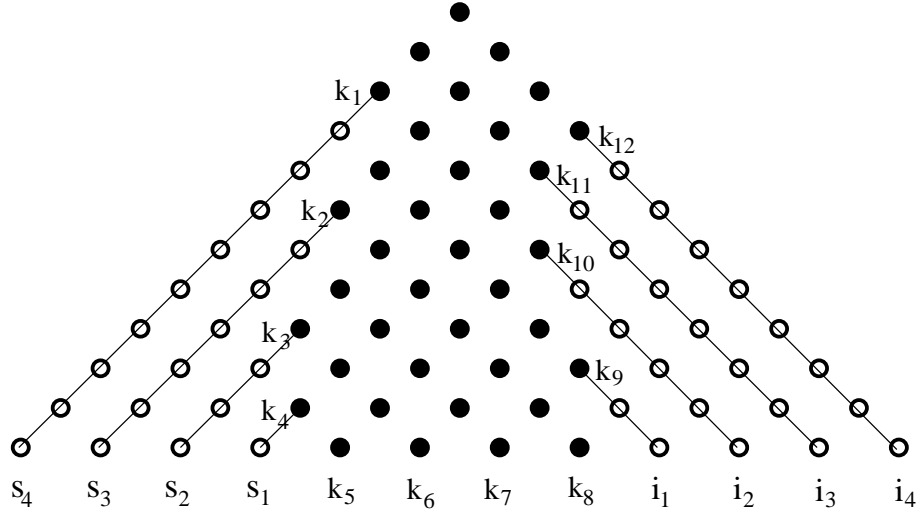}}}
\end{center}
\label{gamma1}
\caption{In this example we have $n=12$, $c=d=4$, $s_1=2, s_2=4,s_3=7,s_4=10$ and $i_1=3,i_2=6,i_3=8,i_4=9$.}   
\end{figure}
    
For $1 \le i_1 \le i_2 \le \ldots \le i_d \le n$ and $1 \le s_1 \le s_2 \le \ldots \le s_c \le n$, 
let $$\gamma(n;k_1,\ldots,k_n;s_1,\ldots,s_c;i_1,\ldots,i_d)$$ (see Figure~\ref{gamma1})
denote the number of partial monotone triangles with $n$ rows that have the
following properties.
\begin{itemize}
\item For $1 \le l \le c$, the first $s_{c+1-l}-1$ entries of the
$l$--th NE--diagonal are missing.
\item The first entry of the
truncated $l$--th NE--diagonal is equal to $k_l$ and does neither have
to be strictly smaller than its neighbour to the right nor weakly smaller than
its SE--neighbour.
\item For $n-d+1 \le l \le n$, the last $i_{l-n+d}-1$ entries in the $l$--th
SE--diagonal are missing.
\item  The last entry of the truncated $l$--th SE--diagonal 
 is equal to $k_l$ and does neither have to be
strictly greater than its neighbour to the left nor weakly greater than its SW--neighbour.
\item The bottom row
(without those entries among $k_1, k_2, \ldots, k_c, k_{n-d+1},
k_{n-d+2}, \ldots, k_n$ that are situated in the bottom row) is
equal to $k_{c+1},k_{c+2},\ldots,k_{n-d}$ and we do not demand
a strict increase for those entries.
\end{itemize}
We want to show that
\begin{multline}
\label{gamma}
 \gamma(n;k_1,\ldots,k_n;s_1,\ldots,s_c;i_1,\ldots,i_d) \\
=
 \Delta^{s_{c}-1}_{k_1} \Delta^{s_{c-1}-1}_{k_2} \ldots
\Delta_{k_c}^{s_1-1} \delta^{i_1-1}_{k_{n-d+1}}
\delta^{i_2-1}_{k_{n-d+2}} \ldots \delta_{k_n}^{i_d-1}
 (-1)^{s_1+\ldots+s_c-c} \alpha(n;k_1,\ldots,k_n)
\end{multline}
by induction with respect to $n$. The case $n=1$ is easy to check.
Otherwise, let $c' \ge 0$ be minimal such that $s_{c'+1} > 1$ and $d'
\ge 0$ be minimal such that $i_{d'+1} > 1$. Then we have the following
recursion for $\gamma$ with respect to $n$.
\begin{multline*}
\gamma(n;k_1,\ldots,k_n;s_1,\ldots,s_c;i_1,\ldots,i_d)  =
\sum_{(l_{c-c'+1},  \ldots,
l_{n-d+d'-1})}^{(k_{c-c'+1},\ldots,k_{n-d+d'})} \\
\gamma(n-1;k_1,\ldots,k_{c-c'},l_{c-c'+1},\ldots,
l_{n-d+d'-1},k_{n-d+d'+1},\ldots,k_n;
\\ s_{c'+1}-1,s_{c'+2}-1,\ldots,s_c-1;i_{d'+1}-1,i_{d'+2}-1,\ldots,i_d-1)
\end{multline*}
By the induction hypothesis, this is equal to
\begin{multline*}
\sum_{(l_{c-c'+1},  \ldots,
l_{n-d+d'-1})}^{(k_{c-c'+1},\ldots,k_{n-d+d'})}
\Delta^{s_{c}-2}_{k_1} \Delta^{s_{c-1}-2}_{k_2} \ldots
\Delta_{k_{c-c'}}^{s_{c'+1} -2} \delta^{i_{d'+1}-2}_{k_{n-d+d'+1}}
\ldots \delta_{k_n}^{i_d-2}
 (-1)^{(s_{c'+1}-1)+\ldots+(s_c-1)-c+c'}
\\
\alpha(n-1;k_1,\ldots,k_{c-c'},l_{c-c'+1},\ldots,
l_{n-d+d'-1},k_{n-d+d'+1},\ldots,k_n).
\end{multline*}
Now, \eqref{delta} implies that this is furthermore equal to
\begin{multline*}
\Delta^{s_{c}-2}_{k_1} \Delta^{s_{c-1}-2}_{k_2} \ldots
\Delta_{k_{c-c'}}^{s_{c'+1} -2} \delta^{i_{d'+1}-2}_{k_{n-d+d'+1}}
\ldots \delta_{k_n}^{i_d-2}
 (-1)^{s_{c'+1}+\ldots+s_c} \\ \Delta_{k_1} \Delta_{k_2} \ldots
 \Delta_{k_{c-c'}} \delta_{k_{n-d+d'+1}} \delta_{k_{n-d+d'+2}} \ldots
 \delta_{k_{n}} (-1)^{c-c'}
\alpha(n;k_1,\ldots,k_n).
\end{multline*}
This concludes the proof of \eqref{gamma}.

\medskip

By the observations from the previous section and \eqref{gamma}, the number of
$(d,n-c)$--trapezoids with $(i_{1},\ldots,i_{d})$ as top row and
whose bottom row consists of the numbers in $\{1,2,\ldots,n\}
\setminus \{s_{1},s_{2},\ldots,s_{c}\}$, arranged in increasing
order, is equal to
\begin{multline*}
 \Delta^{s_{c}-1}_{k_1} \Delta^{s_{c-1}-1}_{k_2} \ldots
\Delta_{k_c}^{s_1-1} \delta^{i_1-1}_{k_{n-d+1}}
\delta^{i_2-1}_{k_{n-d+2}} \ldots \delta_{k_n}^{i_d-1} \\
\left. \times (-1)^{s_1+\ldots+s_c-c} \alpha(n;k_1,\ldots,k_n)
\right|_{(k_{1},\ldots,k_{n})=((c+1)^c,c+1,c+2,\ldots,n-d,(n-d)^d)}.
\end{multline*}
Finally, by the expansion \eqref{expansion}, this number is equal to
\begin{multline*}
\sum_{t_{1}=1}^{n} \ldots \sum_{t_{c}=1}^{n}
\sum_{j_{1}=1}^{n} \ldots
\sum_{j_{d}=1}^{n} A(n;t_{1},\ldots,t_{c};j_{1},\ldots,j_{d})
(-1)^{t_{1}+\cdots+t_{c}+s_{1}+\cdots+s_{c}}
\prod_{l=1}^{c} \Delta_{k_{l}}^{s_{c+1-l}-1} \binom{k_{l}-c-1}{t_{c+1-l}-1}
 \\  \left. \times
\prod_{l=n-d+1}^{n} \delta_{k_{l}}^{i_{l-n+d}-1}
\binom{k_{l}-n+d-2+j_{l-n+d}}{j_{l-n+d}-1} \right|_{(k_{1},\ldots,k_{c},k_{n-d+1},\ldots,k_{n})=((c+1)^{c},(n-d)^{d})}.
\end{multline*}
Since $\Delta_x \binom{x}{n} =
\binom{x}{n-1}$ and $\delta_x \binom{x}{n} = \binom{x-1}{n-1}$, this is equal
to
\begin{multline*}
\sum_{t_{1}=1}^{n} \ldots \sum_{t_{c}=1}^{n}
\sum_{j_{1}=1}^{n} \ldots
\sum_{j_{d}=1}^{n} A(n;t_{1},\ldots,t_{c};j_{1},\ldots,j_{d})
(-1)^{t_{1}+\cdots+t_{c}+s_{1}+\cdots+s_{c}}
\prod_{l=1}^{c}  \binom{k_{l}-c-1}{t_{c+1-l}-s_{c+1-l}}
 \\  \left. \times
\prod_{l=n-d+1}^{n}
\binom{k_{l}-n+d-1+j_{l-n+d}-i_{l-n+d}}{j_{l-n+d}-i_{l-n+d}}
\right|_{(k_{1},\ldots,k_{c},k_{n-d+1},\ldots,k_{n})=((c+1)^{c},(n-d)^{d})} \\
= \sum_{t_{1}=1}^{n} \ldots \sum_{t_{c}=1}^{n}
\sum_{j_{1}=1}^{n} \ldots
\sum_{j_{d}=1}^{n} A(n;t_{1},\ldots,t_{c};j_{1},\ldots,j_{d})
(-1)^{t_{1}+\cdots+t_{c}+s_{1}+\cdots+s_{c}} \\ \times
\prod_{l=1}^{c}  \binom{0}{t_{c+1-l}-s_{c+1-l}}
\prod_{l=n-d+1}^{n}
\binom{j_{l-n+d}-i_{l-n+d}-1}{j_{l-n+d}-i_{l-n+d}}
\end{multline*}
As $\binom{0}{k}=\delta_{0,k}=\binom{k-1}{k}$ for all integers
$k$, this simplifies to
$A(n;s_{1},\ldots,s_{c};i_{1},\ldots,i_{d})$ and concludes the
proof of Theorem~\ref{7}.

\section{A formula for $A_{n,i,j}$}
\label{relate}

In order to prove \eqref{anij} we make use of the identity
$$
\alpha(n;k_1,\ldots,k_n) = (-1)^{n-1}
\alpha(n;k_2,\ldots,k_n,k_1-n),
$$
which can be found in \cite[Lemma~5]{newproof}. It implies
$$
\alpha(n;k_1,\ldots,k_n) = (-1)^{t n - t}
\alpha(n;k_{t+1},\ldots,k_n,k_1-n,k_2-n,\ldots,k_t-n)
$$
for $t \ge 0$. Moreover, for every integer $z$,
$$
\alpha(n;k_{1},\ldots,k_{n}) = \alpha(n;k_{1}+z,\ldots,k_{n}+z).
$$

\medskip

Therefore,
\begin{multline*}
\sum_{s_1=1}^{n} \dots \sum_{s_c=1}^{n} \sum_{i_1=1}^{n} \dots
\sum_{i_d=1}^{n} A(n;s_1,\ldots,s_c;i_1,\ldots,i_d)
(-1)^{s_1+\ldots+s_c+c} \\ \times \prod_{l=1}^{c}
\binom{k_l-c-1}{s_{c+1-l}-1} \prod_{l=n-d+1}^{n}
\binom{k_l-n+d-2+i_{l-n+d}}{i_{l-n+d}-1} \\
= \alpha(n;k_1,\ldots,k_c,c+1,c+2,\ldots,n-d,k_{n-d+1},\ldots,k_n)
\\ = (-1)^{t n -t}
\alpha(n;k_{t+1},\ldots,k_c,c+1,c+2,\ldots,n-d,k_{n-d+1},\ldots,k_n,k_1-n,k_2-n,\ldots,k_t-n)
\\ =
(-1)^{t n-t} \\ \times
\alpha(n;k_{t+1}-t,\ldots,k_c-t,c+1-t,c+2-t,\ldots,n-d-t,k_{n-d+1}-t,\ldots,k_n-t,k_1-n-t,\ldots,k_t-n-t)
\\
= \sum_{s_{1}=1}^{n} \ldots \sum_{s_{c-t}=1}^{n}
\sum_{i_{1}=1}^{n} \ldots \sum_{i_{d+t}=1}^{n}
A(n;s_{1},\ldots,s_{c-t};i_{1},\ldots,i_{d+t})
(-1)^{s_{1}+\ldots+s_{c-t}+c+t n} \\
\times \binom{k_{t+1}-c-1}{s_{c-t}-1}
\binom{k_{t+2}-c-1}{s_{c-t-1}-1} \cdots
\binom{k_{c}-c-1}{s_{1}-1} \\
\times \binom{k_{n-d+1}-n+d-2+i_{1}}{i_{1}-1} \ldots
\binom{k_{n}-n+d-2+i_{d}}{i_{d}-1} \\
\binom{k_{1}-2n+d-2+i_{d+1}}{i_{d+1}-1}
\binom{k_{2}-2n+d-2+i_{d+2}}{i_{d+2}-1} \cdots
\binom{k_{t}-2n+d-2+i_{d+t}}{i_{d+t}-1}.
\end{multline*}
This implies that
\begin{multline*}
\sum_{s_{c-t+1}=1}^{n} \sum_{s_{c-t+2}=1}^{n} \ldots
\sum_{s_c=1}^{n} A(n;s_1,\ldots,s_c;i_1,\ldots,i_d)
(-1)^{s_{c-t+1}+\ldots+s_c} \binom{k_1-c-1}{s_c-1} \cdots
\binom{k_t-c-1}{s_{c-t+1}-1} \\ = \sum_{i_{d+1}=1}^{n}
\sum_{i_{d+2}=1}^{n} \ldots \sum_{i_{d+t}=1}^{n}
A(n;s_1,\ldots,s_{c-t};i_1,\ldots,i_{d+t}) (-1)^{t n} \\
\times \binom{k_{1}-2n+d-2+i_{d+1}}{i_{d+1}-1} \cdots
\binom{k_{t}-2n+d-2+i_{d+t}}{i_{d+t}-1}.
\end{multline*}
By the Chu-Vandermonde summation, the right--hand side is equal
to
\begin{multline*}
\sum_{i_{d+1}=1}^{n} \sum_{i_{d+2}=1}^{n} \ldots
\sum_{i_{d+t}=1}^{n} \sum_{s_c=1}^{i_{d+1}}
\sum_{s_{c-1}=1}^{i_{d+2}} \ldots \sum_{s_{c-t+1}=1}^{i_{d+t}}
A(n;s_1,\ldots,s_{c-t};i_1,\ldots,i_{d+t}) (-1)^{t n} \\
\times \binom{k_1-c-1}{s_c-1}
\binom{-2n+d-1+i_{d+1}+c}{i_{d+1}-s_c} \binom{k_2-c-1}{s_{c-1}-1}
\binom{-2n+d-1+i_{d+2}+c}{i_{d+2}-s_{c-1}} \\ \cdots
\binom{k_t-c-1}{s_{c-t+1}-1} \binom{-2n+d-1+i_{d+t}+c}{i_{d+t} -
s_{c-t+1}}.
\end{multline*}
Consequently,
\begin{multline*}
A(n;s_1,\ldots,s_c;i_1,\ldots,i_d)  =
\sum_{i_{d+1}=s_c}^{n} \sum_{i_{d+2}=s_{c-1}}^{n} \ldots
\sum_{i_{d+t}=s_{c-t+1}}^{n} A(n;s_1,\ldots,s_{c-t};i_1,\ldots,i_{d+t}) \\
\times  (-1)^{s_{c-t+1}+\ldots+s_c+ t n }
\binom{-2n+d-1+i_{d+1}+c}{i_{d+1}-s_c} \cdots
\binom{-2n+d-1+i_{d+t}+c}{i_{d+t} - s_{c-t+1}}.
\end{multline*}
Equivalently, by using $\binom{n}{k} = (-1)^{k} \binom{k-n-1}{k}$,
\begin{multline}
\label{circuit}
A(n;s_1,\ldots,s_c;i_1,\ldots,i_d)  =
\sum_{i_{d+1}=s_c}^{n} \sum_{i_{d+2}=s_{c-1}}^{n} \ldots
\sum_{i_{d+t}=s_{c-t+1}}^{n} A(n;s_1,\ldots,s_{c-t};i_1,\ldots,i_{d+t}) \\
\times  (-1)^{i_{d+1}+\ldots+i_{d+t}+ t n }
\binom{2n-c-d-s_{c}}{i_{d+1}-s_c} \cdots
\binom{2n-c-d-s_{c-t+1}}{i_{d+t} - s_{c-t+1}}.
\end{multline}

In the special case $t=1$, $c=2$ and $d=0$, this gives
\begin{equation}
\label{relation} A(n;s_{1},s_{2};-) = \sum_{i_{1}=s_{2}}^{n}
(-1)^{n+i_{1}} \binom{2n-2-s_{2}}{i_{1}-s_{2}} A(n;s_{1};i_{1}).
\end{equation}
We fix the following notation: $A_n:=A(n;-,-)$,
$A_{n,i}:=A(n;i;-)$, $B_{n,i,j}=A(n;i;j)$ and $A_{n,i,j} =
A(n;i,j;-)$. That is, $A_n$ is the total number of $n \times n$
alternating sign matrices, $A_{n,i}$ is the number of $n \times n$
alternating sign matrices where the unique $1$ in the bottom
(equivalently top) row is in column $i$, $B_{n,i,j}$ is the number
of $n \times n$ alternating sign matrices where the unique $1$ in
the top row is in column $j$ and the unique $1$ in the bottom
row is in column $i$ and, finally, $A_{n,i,j}$ is the number of
monotone triangles with bottom row
$(1,2,\ldots,i-1,i+1,\ldots,j-1,j+1,\ldots,n)$. (Or, equivalently,
$A_{n,i,j}$ is the number of monotone $(2,n)$--trapezoids with top
row $(i,j)$ and bottom row $(1,2,\ldots,n)$. In the introduction,
an interpretation of $A_{n,i,j}$ in terms of alternating sign
matrices is provided.)

\medskip

Stroganov \cite[Formula (34)]{stroganov} has shown that
$$
B_{n,i,j} - B_{n,i-1,j-1} = \frac{1}{A_{n-1}} \left( A_{n-1,i-1}
\left( A_{n,j} - A_{n,j-1} \right) + A_{n-1,j-1} \left( A_{n,i} -
A_{n,i-1} \right) \right).
$$
(See \cite{colomo} for generalizations of this result.)
By the combinatorial interpretation of the numbers $B_{n,i,j}$ we
have $B_{n,1,j}=A_{n-1,j-1}$ and thus
\begin{multline*}
B_{n,i,k} = A_{n-1,k-i} + \sum_{l=2}^{i} \left( B_{n,l,k-i+l} -
B_{n,l-1,k-i+l-1} \right) \\
= A_{n-1,k-i} + \frac{1}{A_{n-1}} \sum_{l=2}^{i} \left(
A_{n-1,l-1} \left( A_{n,k-i+l} - A_{n,k-i+l-1} \right) +
A_{n-1,k-i+l-1} \left( A_{n,l} - A_{n,l-1} \right) \right).
\end{multline*}
If we define $A_{n,l} = 0$ for $l \notin \{1,2,\ldots,n\}$ then
$$
B_{n,i,k} = \frac{1}{A_{n-1}} \sum_{l=1}^{i} \left(
A_{n-1,l-1} (A_{n,k-i+l} - A_{n,k-i+l-1})  +
A_{n-1,k-i+l-1} (A_{n,l} - A_{n,l-1}) \right).
$$
On the other hand, by  \eqref{relation}, we have
\begin{equation}
\label{relation1} A_{n,i,j} = \sum_{k=j}^n (-1)^{n+k}
\binom{2n-2-j}{k-j} B_{n,i,k}.
\end{equation}
This finally implies \eqref{anij}.

\section{Remarks} \label{remarks}

In this paper we have shown that the number
$A(n;s_{1},\ldots,s_{c};i_{1},\ldots,i_{d})$ of monotone $(d,n-c)$--trapezoids
with prescribed top and bottom row (or, equivalently, the number of certain
$(n-c-d,n)$--partial alternating sign matrices) appears as a coefficient of a
specialization of $\alpha(n;k_{1},\ldots,k_{n})$ with respect to a certain
(binomial) polynomial basis. We have used this to provide a formula for
$A(n;s_{1},s_{2};-)$. Formulas for $A(n;s_{1};-)$ and $A(n;s_{1};i_{1})$
were previously known. This raises the question of whether there
exist similar formulas for the numbers $A(n;s_{1},\ldots,s_{c};i_{1},\ldots,i_{d})$
if $c+d \ge 3$.

\medskip

Formula \eqref{circuit} shows that it suffices to restrict
our attention to the coefficients $$A(n;-,i_1,\ldots,i_d).$$ Next
we (re--)derive a system of linear equations for these numbers. (This
computation has already appeared in \cite{fischerromik}.) To this
end we need a further identity for $\alpha(n;k_{1},\ldots,k_{n})$:
observe that there is a one--to--one correspondence between the
number of monotone triangles with bottom row
$(k_{1},\ldots,k_{n})$ and the number of monotone triangles with
bottom row $(-k_{n},-k_{n-1},\ldots,-k_{1})$ and therefore
$$
\alpha(n;k_{1},\ldots,k_{n}) = \alpha(n;-k_{n},-k_{n-1},\ldots,-k_{1}).
$$
Consequently,
\begin{multline*}
\alpha(n;1,2,\ldots,n-d,k_{n-d+1},k_{n-d+2},\ldots,k_{n}) \\
= \alpha(n;-k_{n},-k_{n-1},\ldots,-k_{n-d+1},-n+d,-n+d+1,\ldots,-1) \\
= (-1)^{d n - n}
\alpha(n;-n+d,-n+d+1,\ldots,-1,-k_{n}-n,-k_{n-1}-n,\ldots,-k_{n-d+1}-n) \\
= (-1)^{d n - n }
\alpha(n;1,2,\ldots,n-d,-k_{n}-d+1,-k_{n-1}-d+1,\ldots,-k_{n-d+1}-d+1).
\end{multline*}

\medskip

In terms of the expansion \eqref{expansion}, this means that
\begin{multline*}
\sum_{i_{1}=1}^{n} \ldots \sum_{i_{d}=1}^{n} A(n;-;i_{1},\ldots,i_{d})
\prod_{l=1}^{d} \binom{k_{n-d+l}-n+d-2+i_{l}}{i_{l}-1} \\
= (-1)^{d n - d} \sum_{j_{1}=1}^{n} \ldots \sum_{j_{d}=1}^{n} A(n;-;j_{1},\ldots,j_{d})
\prod_{l=1}^{d} \binom{-k_{n-d+l}-n-1+j_{d+1-l}}{j_{d+1-l}-1}.
\end{multline*}
By the Chu--Vandermonde summation, we have
$$
\binom{-k_{n-d+l}-n-1+j_{d+1-l}}{j_{d+1-l}-1} =
\sum_{i_{l}=1}^{j_{d+1-l}} \binom{k_{n-d+l}-n+d-2+i_{l}}{i_{l}-1} (-1)^{i_{l}-1}
\binom{-2n-1+j_{d+1-l}+d}{j_{d+1-l}-i_{l}}
$$
and, consequently,
$$
A(n;-;i_{1},\ldots,i_{d}) \\
= \sum_{j_{1}=i_{d}}^{n} \ldots \sum_{j_{d}=i_{1}}^{n} (-1)^{d n + i_{1} +
  \ldots + i_{d}} A(n;-;j_{1},\ldots,j_{d}) \prod_{l=1}^{d}
\binom{-2n-1+j_{d+1-l}+ d }{j_{d+1-l} - i_{l}}.
$$
Equivalently, by using $\binom{n}{k} = (-1)^{k} \binom{k-n-1}{k}$,
\begin{equation}
\label{system}
A(n;-;i_{1},\ldots,i_{d})
= \sum_{j_{1}=i_{1}}^{n} \ldots \sum_{j_{d}=i_{d}}^{n} (-1)^{d n + j_{1} +
  \ldots + j_{d}} A(n;-;j_{d},\ldots,j_{1}) \prod_{l=1}^{d}
\binom{2n-i_{l} - d }{j_{l} - i_{l}}
\end{equation}
where $1 \le i_{1}, i_{2}, \ldots, i_{d} \le n$. (This is
Conjecture~5 from \cite{fischerromik}.)

\begin{rem} In fact, we have used the principal that certain identities for $\alpha(n;k_{1},\ldots,k_{n})$
  translate to identities for the coefficients
  $A(n;s_{1},\ldots,s_{c};i_{1},\ldots,i_{d})$. Another (much less
  interesting)
  example in this respect is
\begin{multline*}
\alpha(n;k_{1},\ldots,k_{c},c+1,\ldots,n-d,k_{n-d+1},\ldots,k_{n}) \\
=
\alpha(n;-k_{n},-k_{n-1},\ldots,-k_{n-d+1},-n+d,\ldots,-c-1,-k_{c},\ldots,-k_{1})
\\
= \alpha(n;n+1-k_{n},n+1-k_{n-1},\ldots,n+1-k_{n-d+1},d+1,\ldots,n-c,n+1-k_{c},\ldots,n+1-k_{1}),
\end{multline*}
which translates in a similar way to
$$
A(n;s_{1},\ldots,s_{c};i_{1},\ldots,i_{d}) =
A(n;i_{1},\ldots,i_{d};s_{1},\ldots,s_{c}),
$$
an identity that is of course (almost) obvious from the
combinatorial interpretation if $s_{1} < s_{2} < \ldots < s_{c}$
and $i_{1} < i_{2} < \ldots < i_{d}$.
\end{rem}

For fixed $n$, the system of linear equations \eqref{system} provides a total of $n^{d}$
linear equations for the $n^{d}$ numbers
$A(n;-;i_{1},\ldots,i_{d})$, $1 \le i_1, i_2, \ldots, i_d \le n$. However, computerexperiments show
that there is a certain dependency under these linear equations.

\medskip

For the case $d=1$ (i.e. the numbers $A(n;-;i)=A_{n,i}$ of the refined
alternating sign matrix theorem), it was shown in \cite{newproof}
that \eqref{system} together with $A_{n,1}=A_{n-1}$
and the symmetry $A_{n,i} = A_{n,n+1-i}$ for $1 \le i \le n$
(which easily follows from the combinatorial interpretation of
the numbers $A_{n,i}$) provides   a system of linear equations that determines
the numbers $A_{n,i}$ uniquely. In \cite{fischerromik} we have
conjectured that this extends to the case $d=2$, i.e. the numbers
$A(n;-;i,j) = A_{n,i,j}$. In this case, we have to add $A_{n,i,n}
= A_{n-1,i}$, the near--symmetry
$$
A_{n,i,j} = A_{n,n+1-j,n+1-i}
$$
for all $i,j \in \{1,2,\ldots,n\}$  unless $(i,j) \notin \{
(n-1,1), (n,2)\}$ and $A_{n,n-1,1} - A_{n-1} = A_{n,n,2}$. If
$i<j$ this near--symmetry also easily follows from the
combinatorial interpretation of the numbers $A_{n,i,j}$, however, for $i \ge
j$, we took some efforts in \cite{fischerromik} to derive it. It
would be of interest to see whether this (conjectural) behaviour
extends to the case $d > 2$. By Cramer's rule this would then at
least provide a determinantal expression for the numbers
$A(n;-,i_{1},\ldots,i_{d})$.

\medskip

As mentioned above, monotone $(d,n-c)$--trapezoids correspond to
certain $(n-c-d,n)$--partial alternating sign matrices. The latter
are generalizations of alternating sign matrices with loosened
columnsum restrictions. Clearly, a next natural step would be to consider
even more general partial alternating sign matrices, where we
loosen the columsum restrictions and the rowsum restrictions. In
particular, the quadruply refined enumeration of alternating sign
matrices where we fix top and bottom row as well as the first and
the last column of the alternating sign matrix, would be of major interest.



\bigskip \noindent
\textsc{
\!\!Ilse Fischer\\
Institut f\"ur Mathematik, Universit\"at Klagenfurt \\
9020 Klagenfurt, and \\
Fakult\"at f\"ur Mathematik, Universit\"at Wien \\
1090 Wien, Austria \\
}
\texttt{Ilse.Fischer@univie.ac.at}

\end{document}